\documentclass[10pt]{amsart}
\usepackage{hyperref}
\usepackage[utf8]{inputenc}
\usepackage{graphicx} 		
\usepackage{amsmath}  		
\usepackage{bibentry}
\nobibliography*
\usepackage{epigraph}
\usepackage[super]{nth}
\usepackage{multirow}
\usepackage[table]{xcolor}
\usepackage{amssymb}

\DeclareMathOperator{\Id}{Id}

\begin{document}
\title[Proof unwinding by programming languages techniques]{Perspectives for proof unwinding by programming languages techniques}
\author{Danko Ilik}
\maketitle
\begin{abstract}
  In this chapter, we propose some future directions of work,
  potentially beneficial to Mathematics and its foundations, based on
  the recent import of methodology from the theory of programming
  languages into proof theory. This scientific essay, written for the
  audience of proof theorists as well as the working mathematician, is
  not a survey of the field, but rather a personal view of the author
  who hopes that it may inspire future and fellow researchers.
\end{abstract}
\section{Introduction}
\setlength{\epigraphwidth}{0.8\textwidth}
\epigraph{\textit{We cannot hope to prove that every definition, every
    symbol, every abbreviation that we introduce is free from
    potential ambiguities, that it does not bring about the
    possibility of a contradiction that might not otherwise have been
    present.}}{N. Bourbaki~\cite{Bourbaki1949}}

\epigraph{\textit{There is an error, I can confess now. Some 40 years
    after the paper was published, the logician Robert M. Solovay from
    the University of California sent me a communication pointing out
    the error. I thought: “How could it be?” I started to look at it
    and finally I realized} [$\ldots$]}{John F. Nash Jr.~\cite{nash}}
\epigraph{\textit{Mathematics arises from all sorts of application or
    insights but in the end must always consist of proofs}, [but]
  \textit{although a real proof is not simply a formalized document
    but a sequence of ideas and insights}, [a] \textit{real proof is
    not something just probably correct.}}{Saunders Mac Lane
 ~\cite{maclane}}
What constitutes a \emph{real proof} is a question at the origin of
mathematical logic. In effect, a real proof is one that can be reduced
to the use of only a few accepted `ideal' principles such as the
axioms for a set theory like ZFC.  And yet certain ideal principles
are far from self-evident as Euclid's axiomatic method would require
them to be. Proof theory was conceived by Hilbert with the program to
further ``recognize the non-contradictory character of all the usual
[ideal] mathematical methods without exceptions''. Around 1960, these
concerns were addressed for the theory of arithmetic and analysis in
the so called \emph{modified} Hilbert program using the early models
of \emph{computation}---proof theory was also pivotal for the
development of computer science (Hilbert's Entscheidungsproblem).

Applying mathematical rigor to formal proofs as the object of study
brought an answer to the question of what a real proof is: a formal
proof can be given semantics in terms of Gödel's system T and
Spector's bar recursion, thereby eliminating logic in favor of pure
computation. However, these early models of computation that were used
to provide the answer to the consistency question, although satisfying
in terms of precision, are cumbersome to use in practice.

Firstly, it is far from clear why the old computational
interpretations are the right ones, for it is often hard to
distinguish them (bar recursion) from brute force search
procedures. We would like to understand the computational answer to
the main consistency questions in terms of modern and more finely
grained computing abstractions, such as the ones developed over the
course of the past four decades in the theory of programming
languages---for research on (natural) models of computation surely did
not end with the invention of the Turing machine and recursive
function theory.

Secondly, the cumbersome machinery, although ingenious, makes it
difficult to address the next level of research questions. Once that
we have the answer to what a real proof \emph{is}, we need to know
what constitutes the \emph{essential data} of a proof---curiously,
this question of finding criteria of greatest simplicity for proofs
was already listed as \nth{24} in Hilbert's famous list of open
problems, but being premature was not included among the ones finally
published~\cite{Thiele2003}.

The title of this chapter refers to \emph{proof unwinding}, the
pioneering research program from the 1950's of Georg Kreisel
\cite{Luckhardt1996}, who started to use the computational approach,
not for foundational purposes, but to extract numerical content from
actual mathematical arguments. We aim at the \emph{working
  mathematician}, a term used by Bourbaki~\cite{Bourbaki1949} who
meant by it a researcher with a pragmatic attitude toward
foundations. The time is ripe for a 
leap forward, both in foundations and unwinding applications. The
present chapter has as goal to propose bringing proof unwinding on a
par with the latest computing abstractions from the theory of
programming languages, with the ambition to turn such
streamlined proof theoretic methods into a toolbox readily used by the
working mathematician, rather than the rare specialist in proof theory
as it has been the case up to now.

\section{New Unwinding Toolbox}
\label{subsec:toolbox}

Conducted with the goal:
\begin{quote}
  \textit{``To determine the constructive (recursive) content or the
    constructive equivalent of the non-constructive concepts and
    theorems used in mathematics''}~\cite{Kreisel1958},
\end{quote}
Kreisel's research program applied the proof theory of the day, namely
Hilbert's $\epsilon$-substitution method, Herbrand's theorem, and the
no-counterexample interpretation, combined with then brand new theory
of recursive functions, to extract new bounds and algorithms from
prima facie ineffective proofs.  But, even in the hands of masters,
the early unwinding methodology was apparently difficult to apply, if
one is to judge from the lapses of time in between applications:
Littlewood's theorem by Kreisel in 1951~\cite{Kreisel1951}, Artin's
proof of Hilbert's \nth{17} problem by Kreisel first in 1957 and again
in 1978~\cite{Delzell1996}), the Thue-Siegel-Roth theorem by Kreisel
and Macintyre in 1982 and Luckhardt in 1989
\cite{KreiselM1982,Luckhardt1989}, Van der Waerden's theorem by Girard
in 1987~\cite{Girard1987}. The unwinding methods are so complex that
there are even doubts cast on some of the results by authorities in
proof theory~\cite{Feferman1996}.

However, there is a more recent application of unwinding to functional
analysis in the \emph{proof mining} program of Kohlenbach
\cite{Kohlenbach2005}. This very successful unwinding program has at
its methodological core the classic unwinding approach using
Kolmogorov's double-negation translation (1929) and Gödel's functional
`Dialectica' interpretation (1941).

In parallel, in constructive mathematics, there have been equally
significant results in the program of \emph{constructive analysis}
\cite{bishopbridges} and \emph{constructive algebra}
\cite{mines88,Lombardi}, although these are primarily theory
reconstruction programs and rely little on direct application of proof
theoretic methods to unwind ineffective proofs.

But, both mining and research in constructive mathematics have not
sought to reap the benefits of notable proof theoretic advances
directly inspired by the theory of programming languages. This theory,
a continuation of the work on the early models of computation, has
arrived at highly abstract notions for structuring programs. We shall
now describe the proof theoretic state-of-the-art for three such proof
unwinding techniques.  This new methodology will be referred to as
the \emph{unwinding toolbox}.

\subsection{Computational Side-Effects}
\label{methodology:control}
The first of these methods concerns \emph{computational
  side-effect}. Namely, since the work of Griffin~\cite{Griffin1990},
it has been known that the principle of proof by contradiction can be
interpreted by a programming language mechanism (a computational
side-effect) for \emph{control operators}. Although, in absence of
mathematical axioms additional to the reductio-ad-absurdum principle,
control operators provide not much more than a very elegant way to
obtain Herbrand's theorem (an important very early result on classical
first-order logic from 1930), in presence of additional axioms like
induction or choice, when Herbrand's theorem no longer holds, one
begins to get new results. For instance, by the use of computational
side-effects, in set theory, Krivine has managed to extend Cohen's
forcing method from the usual sets of conditions to realizability
algebras~\cite{Krivine}.

However, the promise that control operators can turn every proof by
contradiction into an effective one is a mirage: there are classically
provable formulas whose effective proof would allow to decide the
Halting problem. Whether an ineffectively proved formula can be
unwound, in general needs to be considered on a case-by-case
basis. Nevertheless, there \emph{are} whole \emph{classes} of formulas
which we know can be unwound upfront, like the class of
$\Pi^0_2$-formulas. \emph{Delimiting} control operators only to
formulas in these classes allows to get a new constructive logic. This
logic still respects the existence property, characteristic of
intuitionistic logic that is at the bases of current constructive
mathematics, but the obtained new constructive logic manages to prove
intuitionistically non-provable principles.

For instance, Herbelin~\cite{HerbelinMP} showed that Markov's
principle (MP), \[\neg\neg\exists x A_0(x)\to \exists x A_0(x),\] where
$x\in\mathbb{N}$ and $A_0$ is quantifier-free, an axiom crucial for
constructive proofs of completeness of first-order logic
\cite{IlikThesis}, can be interpreted with the help of a computational
side-effect known as \emph{(delimited) exceptions}. The author further
showed~\cite{Ilik2010} that the double negation shift principle (DNS),
\[\forall x \neg\neg A(x) \to \neg\neg\forall x A(x),\] where
$x\in\mathbb{N}$, a principle crucial for the interpretation of the
classical axiom of choice---and that has only been interpreted before
by the \emph{generally}-recursive schema of bar recursion---can be
interpreted computationally by a generalization of the exceptions
effect to so called delimited continuations, or \emph{delimited
  control operators}. The key
observation from these results is that---when
delimited---computational side-effects like control operators can be
used to unwind ineffective proofs and at the same time not run into
non-decidability problems. The newly obtained logics are intermediate
logics, in between classical and intuitionistic logic, and take the
best of both worlds.

These results are controversial from the point of view of the orthodox
constructive mathematician who is used to intuitionistic logic as
first codified by Heyting's analysis of Brouwer's work in
intuitionistic mathematics. Namely, the only previous computational
interpretation of MP were either trivial (as given by Gödel's
functional interpretation) or proceeded by unbounded search. As for
DNS, the computational interpretation was only given by the
generally-recursive bar recursion schema, whose termination must be
ensured by Brouwer's bar induction or continuity principle. As
unbounded/general recursion can lead to an inconsistent formal system,
intuitionists have been understandably wary of accepting these
principle. By replacing the mentioned computational interpretations by
computationally meaningful realizers, we not only propose to
intuitionists to reconsider the constructivity of principles like MP
and DNS, but we are re-establishing the link between modern proof
theory and one of the offspring of Hilbert's proof theory, the theory
of programming languages.

A further principle interpreted in this way, in a joint work of Nakata
and the author~\cite{IlikN2014}, was the open induction principle,
\[\forall \alpha(\forall \beta<\alpha(\beta\in U) \to \alpha\in U) \to
\forall \alpha(\alpha\in U),\]
where $\alpha,\beta$ range over Cantor space and $U$ is open. This
principle is the only known equivalent form of the axiom of choice
that is stable under double-negation translation (even if we replace
Cantor space by Baire space). The principle is also interesting for
combinatorics, where it leads to a direct version of Nash-Williams'
proof of Kruskal's tree theorem~\cite{Veldman2001}, as well as in
algebra where it is used to replace Zorn's lemma~\cite{Schuster}.

We finally mention a last result~\cite{Ilik2014}, still in review, on
the nature of the programming-language inspired proof rules. It
concerns higher-type primitive recursion---Gödel's system T---versus
general recursion---Spector's bar recursion. Namely, already in 1979,
Schwichtenberg has shown that bar recursion of type 0 and 1 does not
allow to define functions beyond system T
\cite{Schwichtenberg1979}. Since a previous analysis of Kreisel
\cite{Kreisel1959} shows that these types are sufficient for all
practical purposes (realizing $\Sigma^0_2$-theorems), it follows that
we know for a long time that we should not need bar recursion for the
computational interpretation of ideal proof principles. What we
proposed is how to circumvent bar recursion and generate System T
terms directly, using delimited control operators as an intermediary
step. This also shows that the extensions of system T with
computational side-effects are in fact conservative extensions. In
order to establish this fact, we relied on \emph{partial evaluation},
the second set of techniques of our unwinding toolbox that we explain
in the following subsection.

\subsection{Partial Evaluation and Formalization}
\label{methodology:interpret}
The second programming languages method that we intend to employ for
proof unwinding concerns formalization of proofs in proof assistant
software and, more specifically, the use of formalized
\emph{(type-directed) partial evaluators}.

The topic of partial evaluators came up in our previous research on
constructive versions of completeness theorems
\cite{IlikThesis}. These logical theorems establish the adequacy of a
given formal system to encode actual mathematical proofs. As it turns
out, and thanks to initial work on the link between normalization proofs
and completeness of intuitionistic logic for Kripke models
\cite{CoquandD1997}, the computational content of proofs of
intuitionistic completeness can be expressed by type-directed partial
evaluation algorithms~\cite{Danvy1996}. Having a rich theory of such
algorithms in the theory of programming languages, allowed to cover
cases of constructive completeness proofs that were beyond the
previous state-of-the art in proof theory. More precisely, we now know
how to partially evaluate (i.e. show constructively completeness for)
not only classical logic~\cite{Ilik2010}, intuitionistic logic with
disjunction~\cite{Ilik2013}, but also simultaneous presence of
delimited control operators and higher-type primitive recursion
\cite{IlikTDPE,Ilik2014} (the second citation is still in review).

The development of these logical meta-theorems was conducted formally,
in the Coq and Agda proof assistants. Since the formalized proofs are
constructive, they can be used to \emph{compute a proof
  transformation} for every actual formalized argument. What this
allows is to perform unwinding of actual mathematical proofs more
directly, by pushing the complexity of doing a manual double-negation
transformation (like done in the classic unwinding approach and used,
for instance, by Kohlenbach in his program of proof mining in
analysis) into the realizability model, that is, into the reduction
mechanism of the proof assistant used.

Proof assistants are most well known for their use in the full
formalization of complex proofs, such as the four-color theorem
\cite{Gonthier2008}, the Kepler conjecture~\cite{Hales2012}, or the
Feit-Thompson theorem~\cite{Gonthier2013}. However, as far as proof
unwinding is concerned, one can in general \emph{avoid} needing a
\emph{fully} formalized version of an actual mathematical proof. It
suffices to notice that lemmas that have a computationally irrelevant
form, such as $\Pi^0_1$, can be simply assumed without proof. A more
refined analysis of computational relevance of formulas can be found
in~\cite{SchwichtenbergW2012} where the classes of so called definite
and goal formulas are isolated. This allows to greatly decrease the
burden of formalizing i.e. we are only dealing with \emph{partial
  formalization} which nonetheless contains as much of algorithmic
content as a full formalization.

The important lesson that we learned from partial evaluation is that
proofs need not be interpreted uniformly, by `oracles' such as bar
recursion that work uniformly (for example the realization of DNS by
bar recursion is agnostic of the concrete formula $A$ in the instance
of DNS). Rather, it is possible to specialize (i.e. partially
evaluate) proofs, even if they are highly ineffective, and, when one
in addition uses a proof assistant like Coq, the specialization of the
(partially) formalized theorem can become automatic. This is one of
the principal advantages of our toolbox over the old toolbox built on
Herbrand's theorem, $\epsilon$-substitution, double-negation- and
A-translation, and functional interpretation: while unwinding, the
mathematician can concentrate on the essential parts of a proof rather
than get lost in manual proof transformations.

\subsection{Type Isomorphisms}\label{methodology:typeiso}
The final third method of our unwinding toolbox concerns \emph{type
  isomorphisms}. Mathematically, this notion is the same as the one of
\emph{constructive cardinality} of sets~\cite{mines88}, saying not
only that sets are of the same size, but moreover that they have
indistinguishable structure. In programming languages theory, the
notion allows to generalize the notion of type assigned to a program,
which allows to test more easily if a programs conforms to a formal
specification~\cite{Rittri}.

The link that brings us to the study of type isomorphisms is
Tarski's high-school algebra problem~\cite{burris04}. This basic
question, asking whether the system of eleven arithmetic equation
taught in high-school suffices to derive all the true equations
between \emph{exponential} polynomials, had taken some time to be
answered in mathematical logic. It turned out that the high-school
system is not complete, as shown by a counter-example of Wilkie in
1981~\cite{wilkie00}, a true statement which is not derivable by only
using the eleven equations. Gurevič further showed that the system
cannot be completed by any additional finite set of axioms
\cite{gurevic90}.

Now, by the Curry-Howard correspondence, formulas of intuitionistic
logic can be seen as types (conjunction correspond to products,
disjunction to coproducts, and implication to exponentials) and proofs
of formulas can be seen as computer programs of the corresponding
type. Following the correspondence, one gets a notion of strong
equivalence, or formula isomorphism, from isomorphism of types. A new
correspondence is thus obtained: the language of formulas is the same
as the language of exponential polynomials---and, moreover---formula
isomorphism generalizes equality of exponential polynomials in the
standard model of positive natural numbers, that is
\[
A \cong B \text{ implies } \mathbb{N}^+ \vDash A=B.
\]

The link that one establishes in this way allows to use the rich
theory on exponential polynomials to obtain proof theoretic
results. For instance, Fiore, Di Cosmo, and Balat, showed that the
non-finite-axiomatizability result of Gurevič also hold for the theory
of type isomorphism~\cite{fiore06}. 
Using results of Richardson, Martin, Levitz, Wilkie, Macintyre, Henson
and Gurevič, the PI proved that although not finitely axiomatizable,
type i.e. formula isomorphism is recursively axiomatizable and
moreover decidable~\cite{sumaxioms}. The value of this unexpected
positive result is still somewhat limited because the existence of a
\emph{practical} decision algorithm for type isomorphism is open.

Nevertheless, even if the meta-theory of type isomorphism has
remaining open questions to be resolved, applications to proof theory
are well under way. Recently, the PI proposed a pseudo-normal form of
types~\cite{explog}, inspired by the decomposition of the exponential
function in exponential fields~\cite{hardy}, called the exp-log normal
form, that allows to decompose the axioms of the notoriously non-local
theory of $\beta\eta$-equality for the lambda calculus with coproduct
type. This equality can be seen as the essence of \emph{identity of
  proofs} for intuitionistic propositional logic with disjunction. An
extension to the first-order case has also been proposed in a joint
work with Brock-Nannestad~\cite{highschool}, where the normal form
appears to produce the first arithmetical hierarchy for formulas of
intuitionistic logic that copes with both quantifiers equally well;
the only previously known hierarchy, the one of Burr~\cite{Burr2000},
covers well only the universal quantifier. This has been a long
standing open problem for constructive logic, although for classical
logic an arithmetical hierarchy exists since the 1930s.

\section{Perspectives}

Today, a paradigm change in proof unwinding is possible, thanks to the
notions from contemporary programming languages theory comprising our
\nameref{subsec:toolbox}. These long-evolved techniques provide
proof-theoretic simplifications of the order that makes them more
accessible even to non-specialists in proof theory.

The overall goal of this undertaking would be to exploit the full potential of
the novel toolbox and apply it, beyond logic itself, to proofs of
landmark results in number theory, combinatorics, and homotopy
theory. In parallel, it would be necessary to address the foundations
of unwinding i.e.  tackle long-standing open questions in the
foundations of constructive mathematics such as identity of proofs and
simplified computational interpretations of semi-intuitionistic
principles. We have thus two sets of objectives, work on applications
and work on foundations.

\subsection*{Objective I --- Applications of Proof Unwinding}
\label{obj:unwinding} The first set of objectives concerns
applications to areas that are important for the `working
mathematician', that is, analysis, number theory, and combinatorics,
as well as an application to unwinding incompleteness theorems in
logic. Objective I would be achieved by tackling three more specific
objectives, called \emph{perspectives}: \nameref{task:mining},
\nameref{task:unwinding}, and \nameref{task:incompleteness}.

\subsection*{Objective II --- Foundations of Proof Unwinding}
\label{obj:nextgen} The second set of objectives concerns work on
foundations of constructive mathematics that are both necessary to
guarantee the soundness of applying unwinding and as an update to the
current foundational theories. The two more specific objectives, or
perspectives, to be tackled are: \nameref{task:hott} and
\nameref{task:metatheory}.

\vfill
The immediate effects of the project would be, on the one hand, to show
that our new proof theoretic methods can be used by the working
mathematician to extract numerical bounds and algorithms from prima
facie ineffective proofs in analysis, number theory, combinatorics,
homotopy theory, and logic, and, on the other hand, to update the
current foundational theories of constructive mathematics with the
powerful computing abstractions that computational side-effects,
partial evaluators and type isomorphisms represent.

In the longer term, we can hope to see the streamlined proof unwinding
methodology becoming an important toolbox across mathematics. We can
also expect to see a synergy of the objectives. For instance, not only
would unwinding efforts across mathematics become possible (Objective
I), but, as the new constructive foundations (Objective II) get
adopted in the community working on proof assistant systems, proof
analysis and development would eventually be carried out even more
efficiently with the help of a proof assistant.

The approach to fulfilling the two objectives would be through carrying
out the five perspectives described in this section. We shall explain each
one of the tasks in the context of its proper state-of-the-art,
objectives, methodology, and feasibility.

\subsection{Perspective 1: Unwinding in Analysis Revisited}
\label{task:mining}
Analysis is essentially the only area of mainstream mathematics to
have benefited from direct application of proof unwinding techniques.
In approximation theory, by using proof theory, Kohlenbach and his
collaborators have managed to obtain explicit moduli of uniqueness,
significantly better than previous ones, for best Chebyshev
approximation, as well as to obtain a first effective rate of strong
unicity in the case of best approximation for the $L_1$-norm
\cite{Kohlenbach1993}. How this works is that first logical
meta-theorems are established~\cite{Kohlenbach2005}, which are on one
hand general enough to be applicable as analytic theorems, and on the
other hand specific enough to enter in a class of statements, such as
the $\Pi^0_2$-class of the arithmetical hierarchy, for which we know
by proof theory that explicit functions or existence witnesses can be
extracted. Moreover, such general logical meta-theorems are not only
good for extracting numerical data from concrete proofs, but also for
analyzing whether a known analytic theorem has optimal form. For
instance, in the fixed point theory for functions of contractive type,
one does not only get effective quantitative forms of theorems, but
one can often also relax the compactness assumption for the metric
space.

Why, then, when Kohlenbach's proof mining approach is already
successful, do we propose a perspective on proof unwinding of analysis? There
are two reasons. The first one is methodological: our form of
unwinding has not been applied outside of logic, and proof mining
provides the perfect test bed to make it grow up in the `real
world'. Second, even if we cannot pretend to analytic skills of the
level of the ones present in mining, we do believe that the general
logical meta-theorems can be unwound in a simpler way; this could lead
to better extracted bounds even if we use the exact same analytic
machinery as in mining.

To explain the difference and simplification mentioned, we briefly
explain how the meta-theorems are established right now. The core is
to show that in classical logic, and in presence of additional axioms
for induction and choice, like the weak Kőnig's lemma, one can turn
the $\forall\exists$ quantifier combination from
$\forall x \exists y A_0(x,y)$, where $A_0$ is a quantifier-free
formula, into an explicit recursive function $f$ such that
$\forall x A_0(x, f(x))$.  One can further extend this to formulas
beyond the strict class $\Pi^0_2$ and allow for instance any number of
additional hypotheses of form $\Pi^0_1$. But, to obtain the recursive
functions $f$, which, as explained before in the section
\nameref{methodology:control}, needs an a priori generally recursive
definition schema, one first has to transform by the double-negation
translation all proofs of the original proof system (Peano arithmetic
+ axiom of choice) into proofs of a (semi-)intuitionistic system. This
is a \emph{non-local} transformation of proofs, and in particular the
meaning of formulas can be changed by the transformation (hence the
restriction to the $\Pi^0_2$-class of formulas). Once
a~(semi-)intuitionistic proof is obtained, one can use Gödel's
functional interpretation to obtain a higher-type primitive recursive
function, possibly also needing Spector's generally-recursive schema
of bar recursion. Actually, more redefined versions of the Dialectica
interpretation (monotone and bounded variants) and of bar recursion
are used in practice.

Now, what our approach offers is first to push the technical
complexity of the double negation translation into the realizability
model based on computational side-effects (ex. control
operators). Since these notions have a well-studied operational
semantics, one can perform a more direct reduction of a proof to a
program or a more direct reading off of witnesses (numeric
bounds). With the additional help of a proof assistant like Coq, this
can be further automatized.

In addition, thanks to the reasons already explained in section
\nameref {methodology:control}, our unwinding method makes it likely
that in fact a pure system T witnessing terms can be extracted from
any concrete proof, circumventing bar-recursion-like schemata
altogether.

A third, orthogonal, improvement to the extraction process will be
offered by use of richer data-types for extracted programs and
bounds. Traditionally, one only uses the `negative' function and
product types. Although these can encode `positive' types (for
instance, sum types $\rho+\sigma$ can be encoded by
$(\mathbb{N}\to \rho)\times(\mathbb{N}\to \sigma)$), encoding leads to
an increase of the \emph{degree} of the type. Simpler and more natural
realizers can thus be extracted from disjunction and other inductively
defined positive predicates.

\subsubsection*{Feasibility for Perspective 1} We will need to cope with
semi-intuitionistic principles that we have not treated before,
notably the weak Kőnig's lemma and the independence of premise
schema. For these, we plan to use \nameref{methodology:control}, like
we have done previously for the open induction principle: the fan
theorem, a positive version of the weak Kőnig lemma, is implied by
open induction. At the level of realizers, it will be necessary to use
{\nameref{methodology:typeiso}} to handle extensionality.

The risk for handling the logical part (meta-theorems) is moderate,
hence it is possible to propose this for a subject of a PhD thesis. As
for obtaining better bounds that the ones already obtained in proof
mining, the risk is higher; in fact, it would be a success even if we
manage to obtain the same bounds, since this would mean that our
toolbox is ready to be used in the following, Perspective 2.

\subsection{Perspective 2: Unwinding in Number Theory and Combinatorics} 
\label{task:unwinding}

In this task, we should bring in our \nameref{subsec:toolbox} to bear on
highly non-effective proofs from number theory and combinatorics. The
concrete goals will be to unwind landmark proofs in these areas, but
what we see as equally important is to arrive at a situation where a
sufficiently interesting intersection of proof theory and the
application domain area is recognized. This kind of objective is only
possible through a combination of expertise, and for its carrying out,
it would be appropriate to engage two post-doctoral researchers, one in
each application domain.

In number theory, we would intend to unwind Thue-Siegel-Roth's theorem on
Diophantine approximations. Saying that an algebraic irrational number
has only finitely many exceptionally good rational number
approximations, this $\Sigma_2$ statement has first been tackled upon
by Kreisel and Macintyre~\cite{KreiselM1982} using technology for
obtaining Herbrand terms. However the combinatorial explosion arising
from use of Herbrand's theorem apparently did not allow to obtain
useful bounds on the number of rational approximations, and only
Luckhardt~\cite{Luckhardt1989} managed to limit the growth of Herbrand
terms in order to obtain such a bound. This bound is essentially the
same as the one obtained by Bombieri and van der Poorten
\cite{Bombieri}.

In this case, even more advanced existing technology like Gödel's
functional interpretation has not been applied. We suspect this is the
case because, in order to apply it, one would first need to perform a
double-negation translation of an actual proof of Thue-Siegel-Roth
into a semi-intuitionistic theory---something possible to do in
principle, but given the sophistication of the original proof, its
translation would be an order of magnitude more complex. We propose
thus to treat it directly using our approach with computational
side-effects, i.e. without a preliminary double-negation translation
followed by a functional interpretation. Technically, our approach can
be seen as a version of the so called modified realizability technique
but where the language of realizers is enriched to contain delimited
control operators.

In combinatorics, we would intend to unwind Szemerédi's theorem saying that
every subset of the natural numbers with positive upper density
contains arithmetic progressions of arbitrary length. Conjectured by
Erdős and Turán in 1936, this statement was only proved by Szemerédi
in 1975 by an ingenious and complex combinatorial argument. In 1977,
Furstenberg provided a proof using ergodic theory. The interest in
giving a better proof of this theorem is still ongoing, and
applications include for instance the recent work of Green and Tao on
arbitrary long arithmetic progressions in the prime numbers
\cite{green2008primes}.

We first intend to address an important special case of the theorem,
the van der Waerden theorem, saying that if we use a finite number of
colors to color the natural numbers, then there is at least one
color containing arbitrarily long arithmetic progressions. The
current upper bound for van der Waerden's number $W(k,r)$, where $r$
is the number of colors and $k$ is the requested length of an
arithmetic progression, was obtained via Szemerédi's theorem and is
due to Gowers~\cite{Gowers}. What is intriguing in this subject is
that the upper bounds appear to be heavy overestimates: for
instance, the bound for $W(3,3)$ is of the order of $10^{14616}$,
while the exact value is $27$.

Girard has previously analyzed Furstenberg and Weiss's proof of van
der Waerden's theorem using cut elimination~\cite{Girard1987}. But,
the bound that he arrived at was essentially the same upper bound
obtained by Furstenberg and Weiss~\cite{Macintyre2005}. We could
attack the problem by using our modified realizability based on
computational side-effects and attempt to partially evaluate the
latest available proofs for Szemerédi's and van der Waerden's
theorem---that would avoid the need for having a fully formalized proof
on hand.

\subsubsection*{Feasibility for Perspective 2} Although the \emph{statements} of
the mentioned theorems in number theory and combinatorics are
arithmetical, their \emph{proofs} are not arithmetical. The risk for
the objectives of this task is to cope with the highly non-effective
nature of proofs, as well as their considerable complexity (see
Szemerédi's diagram of lemmas from his proof in
\cite{Szemeredi}). After all, proofs of the corresponding theorems
have brought Fields medals to both Roth and Szemerédi. The main proof
theoretic question is which kind of ideal principles are at the core
of arguments and can we provide a direct constructive justification
for them. Sometimes, as in the case of Kruskal's theorem, another
statement of Ramsey theory, the link to the open induction principle
(analyzed previously in joint work with Nakata~\cite{IlikN2014}) turns
out to be direct~\cite{Veldman2001}.

We intend to use proof assistant technology and partial formalization
to cope with the complexity of proofs. Concerning mathematical risk,
given a choice of motivated post-doctoral researchers to work on this
topics, I would say that it is medium.  Work on ergodic theory done in
the previous \nameref{task:mining} would serve as preparatory work and
would help to further mitigate the risk.  This task demands more
resources than the other ones.

\subsection{Perspective 3: Unwinding Incompleteness Theorems}
\label{task:incompleteness}

A statement is said to be independent from a theory if it can neither
be proved nor disproved from the axioms of the theory. The
\emph{incompleteness} phenomenon is the fact that for \emph{any}
theory, under the assumption that it is consistent, there exist
statements that are independent of the theory. One might wonder what
is the nature of these statements, and whether they are relevant in
practice. Indeed, the first such statement discovered by Gödel in 1931
has an `artificial' flavor since it encodes the Epimenides' liar
paradox. But, later, natural examples from Ramsey theory have been
found, first by Paris and Harrington~\cite{paris1977mathematical}, and
include important results like Kruskal's tree theorem. Finding
concrete mathematical incompleteness statements is nowadays a fruitful
field of research led by Friedman~\cite{Friedman}.

However, what we find especially interesting is the \emph{limit} at
which a statement starts to become independent from a theory. This
phenomenon, called \emph{phase transition} by analogy to
thermodynamics, happens when the \emph{provability} of a theorem,
taking a rational number as parameter, depends on the \emph{value} of
this rational parameter. For instance, a parametrized version of
Kruskal's theorem can be provable in Peano arithmetic (PA) below a
certain value of the parameter, and becomes independent above that
value---this is in fact a real number, often found by use of analytic
combinatorics, and provides a measure of the strength of the axiom
system. Phase transitions are a research program led by Weiermann
\cite{Weiermann2005}.

In this task, we propose to develop a method for unwinding
incompleteness theorems and phase transition phenomena based on
programming language theory. The idea is that an incompleteness
theorem,
$\text{PA}\not\vdash\bot \rightarrow \text{PA}\not\vdash
\text{Con(PA)}$,
saying that no consistent formal system (in this case, Peano
arithmetic (PA)) can prove its own consistency, can be rephrased
positively as
$\text{PA}\vdash \text{Con(PA)} \rightarrow \text{PA}\vdash\bot$.
Translated in programming languages terms,
$\text{PA}\vdash \text{Con(PA)}$ expresses the possibility of writing
an interpreter for Gödel's system T inside system T itself---that is,
a \emph{self-interpreter}. Self-interpreters have not only been
studied in programming languages theory, but they are a standard way
to bootstrap a compiler for a programming language.

Nevertheless, self-interpreters are usually written for a
Turing-complete languages like Scheme and ones without a strong typing
discipline. If one adds a type system on top of Scheme one can
retrieve system T in its $\lambda$-calculus formulation. There are
recent intriguing results on typed self-interpreters. Brown and
Palsberg have recently constructed the first \emph{typed}
self-interpreter~\cite{brown2015self}; their target was Girard's
system U, and this is still `acceptable', since system U is known to
be inconsistent as a logical system. But, their latest result concerns
Girard's system F$_\omega$~\cite{brown2016self}, which is a
higher-order logic and considered to be consistent.

In this task, we would first investigate whether it is possible to
construct a self-interpreter for system T. For the purpose of the
paper~\cite{Ilik2014}, we have already developed a formally verified
interpreter for system T$^+$ inside Martin-Löf type theory. Since this
type theory has a realizability model based on system T, one comes
close to having a self-interpreter.  We would also have to study the
recent results of Brown and Palsberg, and attempt to retrieve their
result for system F$_\omega$ in system T.

\subsubsection*{Feasibility for Perspective 3} The proposed methodology involves a
frontal attack on consistency of PA. Although the risk is high, the
fact that the prior works of Brown and Palsberg, and the PI, all
involve formally checked proofs, gives us some confidence. If our
effort succeeds, the gain one may have will be equally high as the
risk. But, even if it turns out to be impossible to write a typed
self-interpreter for T, we can aim to obtain solid results on
interpreting Weiermann's phase transition, and hence characterizing
the strength of a formal system, in terms of notions that are equally
natural from the point of view of computation as analytic
combinatorics are.

\subsection{Perspective 4: Identity of Proofs and Homotopy Type Theory} 
\label{task:hott}

Formal proofs are combinatorial objects meant to encode a fully
correct mathematical argument, going down to the smallest details, but
that makes it difficult to spot the most essential parts of an
arguments. Curiously, finding ``criteria of simplicity, or proof of
the greatest simplicity of certain proofs'' was already part of
Hilbert's program, who even planned to include it as the \nth{24} in
his famous list of open problems~\cite{Thiele2003}. In particular,
Hilbert asked for a procedure to decide when two given proof are
essentially the same. This problem known as \emph{identity of proofs}
is still open~\cite{Dosen2003}.

In this task, we would start by tackling the identity of proofs for
constructive logic, before proceeding to a vast generalization of it,
the computational interpretation of Voevodsky's univalence axiom in
homotopy type theory~\cite{Coquand2014} in the case when the
underlying definitional equality has been strengthened to decide
identity of proofs i.e. to convertibility modulo isomorphism for
dependent types.

For intuitionistic \emph{propositional} logic, the difficulty of
deciding identity of proofs comes from the simultaneous presence of
disjunction and implication. Nevertheless, if we follow the analogy
between formulas, types, and exponential polynomials, explained in
section \nameref{methodology:typeiso}, we can re-express the problem
precisely as that of the effective decidability of the
$\beta\eta$-equational theory for the lambda calculus with
coproducts. We have recently proposed a first step in this direction
by showing how to decompose the equational theory for terms, by the
use of the exp-log normal form for types in order to enlarge the
$\beta\eta$-congruence classes of terms~\cite{explog} (in review).

This exp-log normal form of types is extensible to the
\emph{first-order} case, when the quantifiers $\forall$ and $\exists$
are also present. Namely, recent work with Brock-Nannestad \cite{highschool}
shows that it leads to an intuitionistic arithmetical hierarchy, a
classification of formulas that was elusive for intuitionistic logic,
even though it has existed for classical logic since the 1930's where
it is at the basis of results like the completeness theorem.

A further question is whether we can make the technique work for
dependent types, an extension of the first-order case. Martin-Löf Type
Theory has dependent types which allow it to have special treatment of
equality. Basic equality between elements $a,a'$ of a type $A$ is
encoded by the identity type for $A$, $\Id_A(a,a')$. Identity of
proofs in this context means extending the notion of definitional
(computational) equality to cope with $\eta$-equality for coproducts
(and other inductive types).

Pursuing generalization even further, we can talk about \emph{identity
  between proofs of identity}, $\Id_{\Id_A}(p,p')$, that, in turn,
endows ever type $A$ with the structure of a groupoid. Iterating this
construction, $\Id_{\Id_{\Id_{\cdots\Id_A}}}$, allows to show that
every type $A$ is in fact endowed with the structure of an
$\infty$-groupoid~\cite{HofmannStreicher}. Using Grothendieck's
correspondence between $\infty$-groupoids and homotopy types has led
Voevodsky to give a homotopy theoretic interpretation of type theory
in his model based on simplicial sets
\cite{kapulkin2012simplicial}. This model satisfies Voevodsky's
\emph{univalence axiom}, generalizing identity of proofs, and
specializing to: equality at the level of propositions, bijection at
the level of sets, categorical equivalence at the level of groupoids,
etc. Adding this axiom on top of Martin-Löf's type theory produces
homotopy type theory, which is a logical system formalizing the
\emph{univalent foundations} of mathematics~\cite{hottbook}.

What we propose to do is to build the convertibility of proof terms
modulo type isomorphism into the definitional equality of Martin-Löf
and homotopy type theory. An identity type then gets to cover equality
between terms of a whole class of isomorphic types instead of only one
type. We hope that in this way it will be possible to strengthen the
notion of \emph{transport of structures} and to show that important
special cases of the univalence axioms satisfy a simple computational
interpretation. The only existing computational interpretation of
homotopy type theory appears in the effective version of the
simplicial set model~\cite{bezem2014model} and works for the standard
(restricted) notion of identity type.

\subsubsection*{Feasibility for Perspective 4} The univalence axiom is known to
imply a form of full functional extensionality in type theory. Given
that extensionality of functions in general is undecidable, the risk
for extending the computational interpretation for the univalence
axiom defined over the notion of identity types strengthened to work
modulo isomorphism is high.  Nevertheless, by strengthening the
underlying definitional equality of the type theory, we hope to
diminish the need for resorting to full functional extensionality and
even address important special cases of univalence more simply than
before. 

As concerns the identity of proofs for the propositional and
first-order case, based on our preliminary investigations of this
area, we would say that the risk is moderate.

\subsection{Perspective 5: A Next Generation of Constructive
  Foundations}
\label{task:metatheory}

This task would serve as an umbrella for more specific but important
problems that need to be tackled in the foundations of constructive
mathematics, as well as an umbrella collecting the foundational
implications of the previous four tasks of this chapter.

For instance, we already know that axioms which are independent of
intuitionistic logic like double negation shift can be safely added to
intuitionistic systems, but we have to establish the outer limits of
the potential given by \nameref{methodology:control}. We need to
develop direct computational interpretations of principles arising
from the work in constructive reverse mathematics, such as the
equivalent forms of the open induction principle~\cite{veldmanarxiv},
the extension of our work~\cite{IlikN2014} to Baire space, and novel
versions of Markov's principle~\cite{fujiwara}.

Another important topic will be to provide a direct constructive proof
of Goodman's theorem. This theorem says that the axiom of choice
presents a conservative extension of higher type Heyting arithmetic
concerning arithmetical formulas; for the meta-theory of constructive
mathematics, it plays the role that Hilbert's $\epsilon$-elimination
theorems play for the proof theory of classical logic. There has
recently been renewed interest about this old result of Goodman by
other researchers as well
\cite{kohlenbach1999note,coquand2013goodman}.

A third important topic will be to find practical decision algorithms
for type isomorphism. As explain in the section
\nameref{methodology:typeiso}, although a decidability result holds
for type isomorphisms~\cite{sumaxioms}, thanks to prior work of
Richardson~\cite{richardson69} and Macintyre~\cite{macintyre81}, it is
not clear at the moment whether a (practical) decision algorithm can
be constructed. Arriving at such an algorithm would not only be useful
for proof theory, but also for symbolic computation.

Finally, we would like to interact with the researchers working on proof
assistant systems like Coq. The logical cores of proof assistants are
lagging behind contemporary proof theory. For instance, program
extraction from proofs in a state-of-the-art proof assistants such as
Coq relies on the simplest possible realizability interpretation, the
so called modified realizability interpretation of
Kreisel. Integrating the techniques from the \nameref{subsec:toolbox}
would be beneficial for users of proof assistants because it would
allow for easier formalization of many apparently ineffective proofs.

\subsubsection{Feasibility for Perspective 5} 
The main challenge for this task is that, when we are interpreting
semi-intuitionistic principles, we are working at the limit of
computability: our realizability models for the classical axiom of
choice refute the internal (formal) version of Church's thesis, but
the external weak Church's rule still holds~\cite{Ilik2014} (in
review). It is thus hard to predict upfront how far the outer limits
of constructive foundations can be extended. As concerns Goodman's
theorem, we think the risk involved is not very high, since after all
this result has been establish by non-direct methods.  Finally, the
risk on finding a practical algorithm deciding type isomorphism is
hard to estimate; but even if we manage to find ones that only work
for special cases, the benefits could spread also beyond proof theory.

\bibliographystyle{alpha}  
\bibliography{KGS-book-chapter-danko.bib}

\newcommand{\etalchar}[1]{$^{#1}$}
\begin{thebibliography}{GAA{\etalchar{+}}13}

\bibitem[BB85]{bishopbridges}
Errett Bishop and Douglas~S. Bridges.
\newblock {\em Constructive Analysis}, volume 279 of {\em Grundlehren der
  mathematischen Wissenschaften}.
\newblock Springer-Verlag Berlin Heidelberg, 1985.

\bibitem[BCH14]{bezem2014model}
Marc Bezem, Thierry Coquand, and Simon Huber.
\newblock A model of type theory in cubical sets.
\newblock In {\em 19th International Conference on Types for Proofs and
  Programs (TYPES 2013)}, volume~26, pages 107--128, 2014.

\bibitem[BNI16]{highschool}
Taus Brock-Nannestad and Danko Ilik.
\newblock An intuitionistic formula hierarchy based on high-school identities.
\newblock {\em arXiv:1601.04876}, 2016.
\newblock Submitted.

\bibitem[Bou49]{Bourbaki1949}
N.~Bourbaki.
\newblock Foundations of mathematics for the working mathematician.
\newblock {\em The Journal of Symbolic Logic}, 14(1):1--8, 1949.

\bibitem[BP]{brown2016self}
Matt Brown and Jens Palsberg.
\newblock Breaking through the normalization barrier: A self-interpreter for
  {F}-omega.
\newblock To appear in Proceedings of the 43rd Annual ACM SIGPLAN-SIGACT
  Symposium on Principles of Programming Languages.

\bibitem[BP15]{brown2015self}
Matt Brown and Jens Palsberg.
\newblock Self-representation in {Girard's} system {U}.
\newblock In {\em Proceedings of the 42nd Annual ACM SIGPLAN-SIGACT Symposium
  on Principles of Programming Languages}, pages 471--484. ACM, 2015.

\bibitem[Bur00]{Burr2000}
Wolfgang Burr.
\newblock Fragments of {Heyting} arithmetic.
\newblock {\em The Journal of Symbolic Logic}, 65(3):1223--1240, 2000.

\bibitem[BVdP88]{Bombieri}
E~Bombieri and AJ~Van~der Poorten.
\newblock Some quantitative results related to roth's theorem.
\newblock {\em Journal of the Australian Mathematical Society (Series A)},
  45(02):233--248, 1988.

\bibitem[BY04]{burris04}
Stanley~N. Burris and Karen~A. Yeats.
\newblock The saga of the high school identities.
\newblock {\em Algebra Universalis}, 52:325--342, 2004.

\bibitem[CD97]{CoquandD1997}
Thierry Coquand and Peter Dybjer.
\newblock Intuitionistic model constructions and normalization proofs.
\newblock {\em Mathematical Structures in Computer Science}, 7(1):75--94, 1997.

\bibitem[Coq13]{coquand2013goodman}
Thierry Coquand.
\newblock About {Goodman}'s theorem.
\newblock {\em Annals of Pure and Applied Logic}, 164(4):437--442, 2013.

\bibitem[Coq14]{Coquand2014}
Thierry Coquand.
\newblock Théorie des types dépendants et axiome d'univalence.
\newblock {\em Séminaire BOURBAKI}, 66(1085), 2014.

\bibitem[Dan96]{Danvy1996}
Olivier Danvy.
\newblock Type-directed partial evaluation.
\newblock In {\em Proceedings of the Twenty-Third Annual ACM SIGPLAN SIGACT
  Symposium on Principles of Programming Languages (POPL'96)}, pages 242--257,
  1996.

\bibitem[Del96]{Delzell1996}
Charles~N. Delzell.
\newblock {K}reisel's unwinding of {A}rtin's proof.
\newblock In {P}iergiorgio {O}difreddi, editor, {\em {K}reiseliana. {A}bout and
  {A}round {G}eorg {K}reisel}, pages 113--246. {A} {K} {P}eters, 1996.

\bibitem[Do{\v{s}}03]{Dosen2003}
Kosta Do{\v{s}}en.
\newblock Identity of proofs based on normalization and generality.
\newblock {\em Bulletin of Symbolic Logic}, 9(4):477--503, 2003.

\bibitem[FCB06]{fiore06}
Marcelo Fiore, Roberto~Di Cosmo, and Vincent Balat.
\newblock Remarks on isomorphisms in typed lambda calculi with empty and sum
  types.
\newblock {\em Annals of Pure and Applied Logic}, 141:35--50, 2006.

\bibitem[Fef96]{Feferman1996}
Solomon Feferman.
\newblock {Kreisel's ``Unwinding'' Program}.
\newblock In {P}iergiorgio {O}difreddi, editor, {\em {K}reiseliana. {A}bout and
  {A}round {G}eorg {K}reisel}, pages 247--273. {A} {K} {P}eters, 1996.

\bibitem[FIN15]{fujiwara}
Makoto Fujiwara, Hajime Ishihara, and Takako Nemoto.
\newblock Some principles weaker than {Markov’s} principle.
\newblock {\em Archive for Mathematical Logic}, 54(7-8):861--870, 2015.

\bibitem[Fri15]{Friedman}
Harvey Friedman.
\newblock {\em Boolean Relation Theory and Incompleteness}.
\newblock Lecture Notes in Logic. ASL Publications, 2015.

\bibitem[GAA{\etalchar{+}}13]{Gonthier2013}
Georges Gonthier, Andrea Asperti, Jeremy Avigad, Yves Bertot, Cyril Cohen,
  Fran{\c c}ois Garillot, St{\'e}phane Le~Roux, Assia Mahboubi, Russell
  O'Connor, Sidi Ould~Biha, Ioana Pasca, Laurence Rideau, Alexey Solovyev,
  Enrico Tassi, and Laurent Th{\'e}ry.
\newblock {A Machine-Checked Proof of the Odd Order Theorem}.
\newblock In Sandrine Blazy, Christine Paulin, and David Pichardie, editors,
  {\em {ITP 2013, 4th Conference on Interactive Theorem Proving}}, volume 7998
  of {\em LNCS}, pages 163--179, Rennes, France, July 2013. {Springer}.

\bibitem[Gir87]{Girard1987}
Jean-Yves Girard.
\newblock {\em Proof theory and logical complexity}, volume~1.
\newblock Bibliopolis, Naples, 1987.

\bibitem[Gon08]{Gonthier2008}
Georges Gonthier.
\newblock Formal proof---the four-color theorem.
\newblock {\em Notices of the AMS}, 55(11):1382--1393, December 2008.

\bibitem[Gow01]{Gowers}
W.T. Gowers.
\newblock A new proof of {Szemerédi}'s theorem.
\newblock {\em Geometric \& Functional Analysis GAFA}, 11(3):465--588, 2001.

\bibitem[Gri90]{Griffin1990}
Timothy~G. Griffin.
\newblock A formula-as-types notion of control.
\newblock In {\em Conf. Record 17th Annual ACM Symp. on Principles of
  Programming Languages, POPL'90, San Francisco, CA, USA, 17-19 Jan 1990},
  pages 47--58, 1990.

\bibitem[GT08]{green2008primes}
Ben Green and Terence Tao.
\newblock The primes contain arbitrarily long arithmetic progressions.
\newblock {\em Annals of Mathematics}, pages 481--547, 2008.

\bibitem[Gur90]{gurevic90}
R.~H. Gurevi\v{c}.
\newblock Equational theory of positive numbers with exponentiation is not
  finitely axiomatizable.
\newblock {\em Annals of Pure and Applied Logic}, 49:1--30, 1990.

\bibitem[Hal12]{Hales2012}
Thomas Hales.
\newblock {\em Dense sphere packings: A blueprint for formal proofs}, volume
  400 of {\em London Mathematical Society Lecture Note Series}.
\newblock Cambridge University Press, 2012.

\bibitem[Har10]{hardy}
Godfrey~Harold Hardy.
\newblock {\em {Orders of Infinity. The `Infinitärcalcül' of Paul Du
  Bois-Reymond}}.
\newblock Cambridge Tracts in Mathematic and Mathematical Physics. Cambridge
  University Press, 1910.

\bibitem[Her10]{HerbelinMP}
Hugo Herbelin.
\newblock An intuitionistic logic that proves {M}arkov's principle.
\newblock In {\em Proceedings of the 25th Annual IEEE Symposium on Logic in
  Computer Science, LICS 2010, 11-14 July 2010, Edinburgh, United Kingdom},
  pages 50--56. IEEE Computer Society, 2010.

\bibitem[HS98]{HofmannStreicher}
Martin Hofmann and Thomas Streicher.
\newblock The groupoid interpretation of type theory.
\newblock In {\em Twenty-five years of constructive type theory ({V}enice,
  1995)}, volume~36 of {\em Oxford Logic Guides}, pages 83--111. Oxford Univ.
  Press, New York, 1998.

\bibitem[Ili10]{IlikThesis}
Danko Ilik.
\newblock {\em Constructive Completeness Proofs and Delimited Control}.
\newblock PhD thesis, École Polytechnique, Palaiseau, France, October 2010.

\bibitem[Ili12]{Ilik2010}
Danko Ilik.
\newblock Delimited control operators prove double-negation shift.
\newblock {\em Annals of Pure and Applied Logic}, 163(11):1549 -- 1559, 2012.

\bibitem[Ili13a]{Ilik2013}
Danko Ilik.
\newblock Continuation-passing style models complete for intuitionistic logic.
\newblock {\em Annals of Pure and Applied Logic}, 164(6):651 -- 662, 2013.

\bibitem[Ili13b]{IlikTDPE}
Danko Ilik.
\newblock Type directed partial evaluation for level-1 shift and reset.
\newblock In Ugo de'Liguoro and Alexis Saurin, editors, {\em {\rm Proceedings
  First Workshop on} Control Operators and their Semantics, {\rm Eindhoven, The
  Netherlands, June 24-25, 2013 }}, volume 127 of {\em Electronic Proceedings
  in Theoretical Computer Science}, pages 86--100. Open Publishing Association,
  2013.

\bibitem[Ili14a]{Ilik2014}
Danko Ilik.
\newblock {An interpretation of the Sigma-2 fragment of classical Analysis in
  System T}.
\newblock {\em arXiv:1301.5089}, 2014.
\newblock Submitted.

\bibitem[Ili14b]{sumaxioms}
Danko Ilik.
\newblock Axioms and decidability for type isomorphism in the presence of sums.
\newblock In {\em Proceedings of the Joint Meeting of the Twenty-Third EACSL
  Annual Conference on Computer Science Logic (CSL) and the Twenty-Ninth Annual
  ACM/IEEE Symposium on Logic in Computer Science (LICS)}, pages 53:1--53:7.
  ACM, 2014.

\bibitem[Ili15]{explog}
Danko Ilik.
\newblock On the exp-log normal form of types.
\newblock {\em arXiv:1502.04634}, 2015.
\newblock Submitted.

\bibitem[IN14]{IlikN2014}
Danko Ilik and Keiko Nakata.
\newblock A direct version of {V}eldman's proof of open induction on {C}antor
  space via delimited control operators.
\newblock {\em Leibniz International Proceedings in Informatics (LIPIcs)},
  26:188--201, 2014.

\bibitem[KLV12]{kapulkin2012simplicial}
Chris Kapulkin, Peter~LeFanu Lumsdaine, and Vladimir Voevodsky.
\newblock The simplicial model of univalent foundations.
\newblock {\em arXiv preprint arXiv:1211.2851}, 2012.

\bibitem[KM82]{KreiselM1982}
Georg Kreisel and Angus MacIntyre.
\newblock Constructive logic versus algebraization {I}.
\newblock In A.S. Troelstra and D.~van Dalen, editors, {\em The L.E.J. Brouwer
  Centenary Symposium}, pages 217--260. North-Holland Publishing Company, 1982.

\bibitem[Koh93]{Kohlenbach1993}
Ulrich Kohlenbach.
\newblock New effective moduli of uniqueness and uniform a priori estimates for
  constants of strong unicity by logical analysis of known proofs in best
  approximation theory.
\newblock {\em Numerical Functional Analysis and Optimization},
  14(5-6):581--606, 1993.

\bibitem[Koh99]{kohlenbach1999note}
Ulrich Kohlenbach.
\newblock A note on {Goodman}'s theorem.
\newblock {\em Studia Logica}, 63(1):1--5, 1999.

\bibitem[Koh05]{Kohlenbach2005}
Ulrich Kohlenbach.
\newblock Some logical metatheorems with applications in functional analysis.
\newblock {\em Transactions of the American Mathematical Society},
  357(1):89--128, 2005.

\bibitem[Kre51]{Kreisel1951}
Georg Kreisel.
\newblock On the interpretation of non-finitist proofs—{Part I}.
\newblock {\em The Journal of Symbolic Logic}, 16(04):241--267, 1951.

\bibitem[Kre58]{Kreisel1958}
Georg Kreisel.
\newblock Mathematical significance of consistency proofs.
\newblock {\em The Journal of Symbolic Logic}, 23(02):155--182, 1958.

\bibitem[Kre59]{Kreisel1959}
Georg Kreisel.
\newblock Interpretation of analysis by means of constructive functionals of
  finite types.
\newblock In Arend Heyting, editor, {\em Constructivity in Mathematics,
  Proceedings of the colloqium held at Amsterdam, 1957}, Studies in Logic and
  The Foundations of Mathematics, pages 101--127. North-Holland Publishing
  Company Amsterdam, 1959.

\bibitem[Kri]{Krivine}
Jean-Louis Krivine.
\newblock On the structure of classical realizability models of {ZF}.
\newblock To appear.

\bibitem[LQ11]{Lombardi}
Henri Lombardi and Claude Quitté.
\newblock {\em Algèbre commutative -- Méthodes constructives}.
\newblock Calvage \& Mounet, Paris, 2011.

\bibitem[Luc89]{Luckhardt1989}
Horst Luckhardt.
\newblock {Herbrand-Analysen zweier Beweise des Satzes von Roth: Polynomiale
  Anzahlschranken}.
\newblock {\em The Journal of Symbolic Logic}, 54(01):234--263, 1989.

\bibitem[Luc96]{Luckhardt1996}
Horst Luckhardt.
\newblock {B}ounds {E}xtracted by {K}reisel {F}rom {I}neffective {P}roofs.
\newblock In {P}iergiorgio {O}difreddi, editor, {\em {K}reiseliana. {A}bout and
  {A}round {G}eorg {K}reisel}, pages 289--300. {A} {K} {P}eters, 1996.

\bibitem[Mac81]{macintyre81}
Angus Macintyre.
\newblock {\em Model Theory and Arithmetic}, volume 890 of {\em Lecture Notes
  in Mathematics}, chapter The laws of exponentiation, pages 185--197.
\newblock Springer Berlin Heidelberg, 1981.

\bibitem[Mac05]{Macintyre2005}
Angus Macintyre.
\newblock The mathematical significance of proof theory.
\newblock {\em Philosophical Transactions of the Royal Society of London A:
  Mathematical, Physical and Engineering Sciences}, 363(1835):2419--2435, 2005.

\bibitem[ML97]{maclane}
Saunders Mac~Lane.
\newblock Despite physicists, proof is essential in mathematics.
\newblock {\em Synthese}, 111(2):147--154, 1997.

\bibitem[MR88]{mines88}
Ray Mines and Fred Richman.
\newblock {\em A course in constructive algebra}.
\newblock Springer, 1988.

\bibitem[PH77]{paris1977mathematical}
Jeff Paris and Leo Harrington.
\newblock A mathematical incompleteness in {Peano} arithmetic.
\newblock {\em Handbook of mathematical logic}, 90:1133--1142, 1977.

\bibitem[Ric69]{richardson69}
Daniel Richardson.
\newblock Solution of the identity problem for integral exponential functions.
\newblock {\em Zeitschrift für mathematische Logik und Grundlagen der
  Mathematik}, 15:333--340, 1969.

\bibitem[Rit91]{Rittri}
Mikael Rittri.
\newblock Using types as search keys in function libraries.
\newblock {\em Journal of Functional Programming}, 1:71--89, 1991.

\bibitem[RS15]{nash}
Martin Raussen and Christian Skau.
\newblock Interview with {Abel} laureate {John F. Nash Jr.}
\newblock {\em European Mathematical Society. Newsletter}, 97:26--31, September
  2015.

\bibitem[Sch79]{Schwichtenberg1979}
Helmut Schwichtenberg.
\newblock On bar recursion of types 0 and 1.
\newblock {\em The Journal of Symbolic Logic}, 44(3), 1979.

\bibitem[Sch13]{Schuster}
Peter Schuster.
\newblock Induction in algebra: A first case study.
\newblock {\em Logical Methods in Computer Science}, 9(3):1--19, 2013.

\bibitem[SW12]{SchwichtenbergW2012}
Helmut Schwichtenberg and Stanley~S. Wainer.
\newblock {\em Proofs and Computations}.
\newblock Perspectives in Logic. Cambridge University Press, 2012.

\bibitem[Sze75]{Szemeredi}
Endre Szemer{\'e}di.
\newblock On sets of integers containing no k elements in arithmetic
  progression.
\newblock {\em Acta Arith}, 27(199-245):2, 1975.

\bibitem[Thi03]{Thiele2003}
Rüdinger Thiele.
\newblock Hilbert's twenty-fourth problem.
\newblock {\em American Mathematical Monthly}, 2003.

\bibitem[{Uni}13]{hottbook}
The {Univalent Foundations Program}.
\newblock {\em Homotopy Type Theory: Univalent Foundations of Mathematics}.
\newblock \url{http://homotopytypetheory.org/book}, Institute for Advanced
  Study, 2013.

\bibitem[Vel01]{Veldman2001}
Wim Veldman.
\newblock An intuitionistic proof of {Kruskal}'s theorem.
\newblock {\em Archive for Mathematical Logic}, 43:215--264, 2001.

\bibitem[Vel14]{veldmanarxiv}
Wim Veldman.
\newblock The principle of open induction on {Cantor} space and the
  approximate-fan theorem.
\newblock {\em arXiv preprint 1408.2493}, 2014.

\bibitem[Wei05]{Weiermann2005}
Andreas Weiermann.
\newblock Analytic combinatorics, proof-theoretic ordinals, and phase
  transitions for independence results.
\newblock {\em Annals of Pure and Applied Logic}, 136, 2005.

\bibitem[Wil00]{wilkie00}
Alex Wilkie.
\newblock On exponentiation -- a solution to {T}arski's high school algebra
  problem.
\newblock {\em Quaderni di Matematica}, 6, 2000.

\end{thebibliography}

\end{document}